\documentclass[12pt, a4paper,reqno]{amsart}

\usepackage{amsmath,amssymb,amsthm,graphicx,amsfonts,enumerate, mathrsfs}
\usepackage[T1]{fontenc}
\usepackage[utf8]{inputenc}

\usepackage{xcolor}
\usepackage{calc}
\definecolor{bleu_sombre}{rgb}{0,0,0.6}  \definecolor{rouge_sombre}{rgb}{0.8,0,0}\definecolor{vert_sombre}{rgb}{0,0.6,0}
\usepackage[plainpages=false,colorlinks,linkcolor=bleu_sombre,
citecolor=rouge_sombre,urlcolor=vert_sombre,breaklinks]{hyperref}

\usepackage{subcaption}

\usepackage[english]{babel}
\usepackage{geometry}

\theoremstyle{plain}
\newtheorem{theorem}{{Theorem}}[section]
\newtheorem*{theorem*}{{Theorem}}
\newtheorem{proposition}[theorem]{Proposition}

\newtheorem*{proposition*}{Proposition}
\newtheorem{corollary}[theorem]{Corollary}
\newtheorem*{corollary*}{Corollary}
\newtheorem{lemma}[theorem]{Lemma}
\newtheorem{assumption}[theorem]{Assumption}
\newtheorem*{lemma*}{Lemma}

\theoremstyle{definition}

\newtheorem*{definition*}{Definition}

\theoremstyle{remark}
\newtheorem{remark}[theorem]{Remark}

\makeatletter

\@addtoreset{equation}{section}
\makeatother

\renewcommand{\leq}{\leqslant}	\renewcommand{\geq}{\geqslant}

\newcommand{\R}{\mathbb{R}}

\newcommand{\dd}{\mathrm{d}}
\renewcommand{\Re}{\mathrm{Re}\,}

\title{Tunneling effect in some singular wells}

\author[Y. Guedes Bonthonneau]{Yannick Guedes Bonthonneau}
\address[Y. Guedes Bonthonneau]{DMA, ENS. 45 rue d’Ulm, Paris.}
\email{yannick.guedes.bonthonneau@ens.fr}

\author[L. Morin]{Léo Morin}
\address[L. Morin]{Laboratoire de Mathématiques Jean Leray (UMR CNRS 6629), Nantes Université, 2 rue de la Houssinière, 44322 Nantes Cedex 03, France}
\email{leo.morin@univ-nantes.fr}

\author[N. Raymond]{Nicolas Raymond}
\address[N. Raymond]{Univ Angers, CNRS, LAREMA, Institut Universitaire de France, F-49000 Angers, France}
\email{nicolas.raymond@univ-angers.fr}

\begin{document}

	\begin{abstract}
	The Schrödinger operator in one dimension is considered. When the electric potential has a linear/$V$-shaped well, an approximation of the groundstate is established at a fixed distance from the well. In this situation, the classical WKB analysis and the localization techniques \emph{à la} Agmon near a non-degenerate minimum do not work. We propose an alternative to the WKB analysis better suited to the shape of the well and we prove that our Ansatz is an exponentially good approximation by means of revisited (and simplified) Agmon estimates. 
	\end{abstract}
	
		\maketitle

	
\section{Introduction and motivation}
	
	\subsection{About semiclassical tunneling in one dimension}
	
This article is devoted to the spectral analysis of Schrödinger operators in one dimension, \emph{i.e.}, operators of the form 
\[ \mathscr{L} = -h^2 \partial_x^2 + V(x)\,,\]
acting on $L^2(I)$, for some interval $I$ and some real-valued, confining, continuous, electric potential $V$. By confining, we mean that, if $I$ is unbounded, we have $\inf V<\liminf_{|x|\to+\infty}V(x)$. We are especially interested in the approximation of the groundstate eigenfunction in the semiclassical limit $h\to 0$, which is equivalent to the large electric field limit. 

Such a study is motivated by the \emph{tunneling effect}, and probably began with the celebrated article by Harrell \cite{Harrell}. This work was devoted to symmetric double-well potentials, and he established formulas describing the first spectral gap in terms of the eigenfunctions, see \cite[Theorem 2.1]{Harrell}. The exponentially small size of this gap is directly related to the exponential smallness of one-well eigenfunctions in the so-called \emph{classically forbidden zone}, where the electric potential is larger than the eigenvalue, see \cite[Theorems 2.8 \& 2.22]{Harrell}. These considerations are then applied to the case of an explicit smooth non-degenerate double-well, where the one-well eigenfunctions are estimated at a fixed distance from the well, see \cite[Section III]{Harrell}. However, for more singular wells, such estimates are not as easy to obtain.

This relation between exponential decay of the one-well eigenfunctions and the spectral gaps was then used by Simon \cite{Simon84a, Simon84b} and Helffer-Sjöstrand \cite{HS, HS85b} in order to consider quantum systems of higher dimensions. The reader might also want to consider the monograph \cite[Chapters 3 \& 6]{DS}, where WKB constructions are performed (in any dimension) in smooth non-degenerate wells and proved to be exponentially good approximations in Appendix A, see also \cite{Helffer} and \cite{Robert}.

In this article we analyze these questions for one-dimensional singular wells. In fact this analysis is also the occasion to present an alternative to the standard WKB constructions, together with a revisited (and simplified) approach to Agmon estimates. In particular, the techniques we present below could be adapted to non-degenerate wells as well, providing slightly better decay estimates.

	\subsection{The problem under consideration}
	
	The purpose of the present article is to approximate the groundstate eigenfunction of the Dirichlet realization of $\mathscr{L}$ when $I=\R_+$, with a smooth confining $V$ attaining a unique minimum at $x=0$ and such that $V'(0)>0$. This simple model appears in the study of the so-called Stark effect, see, for instance, \cite{CKPRS22}. It is also motivated by the analysis of broken waveguides (or asymptotically flat triangles), see \cite{DR12, OB15} where effective operators involve the Dirichlet realization of $\mathscr{L}$ on a half-line. Hopefully, our results will help establish accurate tunneling estimates by means of dimensional reductions in these contexts, see, for instance, \cite[Section 1.6]{OB15}.  
	
Let us gather our main assumptions.
\begin{assumption}\label{hyp.0}
We consider a non-negative smooth function $V$ on $[0,+\infty)$ such that $V(0)=0$ and $V'(0)>0$. We assume that it has a unique minimum at $0$ and that it is not attained at infinity. For $\lambda>0$ small enough, we denote by $x_\lambda>0$ the unique solution to $V(x)=\lambda$.  	
\end{assumption}

Under Assumption \ref{hyp.0}, we know that $\mathscr{L}$ has discrete eigenvalues below the essential spectrum, which we denote by
\[
0< \lambda_1 < \lambda_2 < \cdots < \lambda_j < \cdots
\]
As $h\to 0$, there are increasingly many such eigenvalues, even if for fixed $h$ there may be only a finite number of them. We are especially interested in the asymptotic behaviour of the normalized groundstate $f$, in the semiclassical limit $h \to 0$. Our main result is the following.

\begin{theorem}\label{thm.linear.wells}
Under Assumption \ref{hyp.0}, for all $j \geq 1$, we have
\begin{equation}\label{eq.lambdajh}
\lambda_j = |z_j|V'(0)^{\frac 23} h^{\frac 23} + \mathcal{O}(h^{\frac 43}),
\end{equation}
where $z_j$ is the $j$-th zero of the Airy function. For all $x >0$, the $L^2$-normalized positive groundstate satisfies, 	in the limit $h \to 0$,
\[
f(x) = \frac{V'(0)^{\frac13}}{2\sqrt{\pi} |\mathrm{Ai}'(z_1)|}\frac{ e^{- \frac 1h S(x)}}{ V(x)^{\frac 14} h^{\frac 16}} \big( 1+ o(1) \big), \qquad S(x) = \int_0^x \sqrt{(V(s)-\lambda_1)_+} \dd s.
\]
\end{theorem}


	\begin{remark}
		Our proof actually shows that the $\lambda_j$ have asymptotic expansions in powers of $h^{\frac23}$. The main novelty in Theorem \ref{thm.linear.wells} is the asymptotic estimate of the groundstate eigenfunction at a fixed distance from its linear minimum. It relies on two independent steps, which are based on classical ideas when considering tunneling effect, namely:
		\begin{enumerate}[--]
			\item construction of exponentially good quasimodes,
			\item sharp exponential decay estimates of the true eigenfunctions and their approximations.
		\end{enumerate}
		The combination of these two ideas allows one to replace $f$ by the quasimode, for which we have the pointwise asymptotic formula. This strategy is now standard, and was used in particular in the setting of non-degenerate wells, as described in \cite{DS}, but also in the case of degenerate wells ($V(x) \sim x^\alpha$ with $\alpha \geq 2$, see \cite{MR88}). In these situations, the exponential decay is given by the \emph{Agmon estimates} \cite{Agmon}, and quasimodes are in the WKB form.
		
		We stress that in our situation, with linear wells, \emph{both the standard WKB analysis and the standard Agmon strategy fail} due to the linear behavior. The strategy of \cite[Théorème 3.1]{MR88} consisting in gluing an Ansatz living at the bottom of the well, \emph{i.e.}, in $\{V(x)<\lambda\}$, and a WKB Ansatz away from the well, does not work in our context due to a more complicated transition of scales. Indeed, the localisation near the bottom of the well is much stronger in the case of linear wells: the characteristic scale near the minimum is $h^{\frac23}\ll h^{\frac12}$ making the connection between the classically allowed and forbidden regions, $\{V(x)<\lambda\}$ and $\{V(x)>\lambda\}$ respectively, more abrupt than for non-degenerate wells. Therefore, we have to provide substitutes to the WKB quasimodes to deal with this transition, and to revisit the usual Agmon estimates.
		
		This yields a robust approach, that could be used in a more general setting. In fact, our contribution gives an alternative to the WKB constructions in the case of non-degenerate wells, and it simplifies the usual Agmon estimates in this case also. For our linear wells, one will use an Airy profile to connect the classically allowed and forbidden regions. Our quasimode construction is reminiscent of Liouville's technique \cite[Section III]{Liouville}, see also \cite[Section 2]{E60} and \cite[Section 4]{O63}. More precisely, our first Ansatz (see \eqref{eq.fNapp} with $N=0$) is Langer's uniform approximation near a turning point \cite{Langer31}, see also \cite[Chapter 11]{Olver74}. Our contribution lies in its iterative improvement and in the Agmon estimates which allow us to compare it with the true groundstate.
		\end{remark}
	
	\subsection{Organization of the article}
This article is organized as follows. Section \ref{sec.2} is devoted to the exposition of the main insights behind the proof of Theorem \ref{thm.linear.wells}. We state the optimal Agmon estimate and introduce the quasimode. In Section \ref{sec.3}, we provide the reader with the detailed proofs of our statements. This section contains the proof of the Agmon estimates, a linear approximation inside the well, the detailed iterative quasimode construction, and the conclusion of the proof of Theorem \ref{thm.linear.wells}. In Section \ref{sec.4}, an application to a tunneling result is discussed. 

	\section{Main insights}\label{sec.2}
	\subsection{Agmon estimates}
	
	A first key ingredient to prove Theorem \ref{thm.linear.wells} is given by the celebrated Agmon estimates (see \cite{Agmon}), which provide exponential decay estimates on the eigenfunctions of Schrödinger operators, of the form
	\begin{equation}\label{eq.Agmon0}
		\| e^{\Phi /h} f \| \leq C h^{-\kappa} \| f \|, \qquad h \in (0,1].
	\end{equation}
	Two natural questions arise: for which weights $\Phi$ can we obtain such estimates, and with which power $\kappa$? This second question, which may seem of lower interest, will actually have implications if one wants to consider potentials of lower regularity.
	
	In our specific situation, we obtain the following Agmon-type elliptic estimate, with optimal weight $\Phi = S$, and a relatively small power $\kappa$. The proof is given in Section \ref{sec.improvedAgmon} (where a slightly stronger version -- which controls 2 derivatives -- is established).
	
	\begin{proposition}\label{prop.Ag}
Let $K$, $E >0$. There exist $h_0$ and $C >0$ such that for all $h \in (0,h_0)$, all $0\leq \lambda \leq Eh^{\frac 23}$ and all $u \in H^1_0(\R_+) \cap H^2(\R_+)$,
\[
\| e^{S / h} u \|_{L^2(0,K)} \leq C h^{- \frac 43} \| e^{S /h} (\mathscr{L}-\lambda)u \| + C h^{- \frac 13} \| u \|, 
\]
where $S(x) = \int_0^x \sqrt{(V(t)-\lambda)_+} \dd t$.
	\end{proposition}
	
	\begin{remark}
	When $u$ is an eigenfunction, we obtain an estimate of the form \eqref{eq.Agmon0} on compact sets, with $\kappa = \frac 13$. The weight $\Phi = S$ is in fact optimal in the sense that it corresponds to the decay announced in Theorem \ref{thm.linear.wells}.
	\end{remark}
	
	The idea behind the Agmon estimates is to use the ellipticity of the conjugate operator $\mathscr{L}_\Phi-\lambda = e^{\frac \Phi h} (\mathscr{L}-\lambda) e^{- \frac \Phi h}$ away from a small neighborhood of $x=0$, for a suitable weight $\Phi$. A straightforward computation gives that, in the sense of quadratic forms,
	\begin{equation}\label{eq.ReLPhi}
	\Re \mathscr{L}_\Phi-\lambda=h^2D_x^2+V(x)-\Phi'(x)^2-\lambda\,,
	\end{equation}
	which naturally leads us to consider the weight 
	\begin{equation}\label{eq.subAgmon}
	\Phi(x)= \sqrt{1-\varepsilon} \left|\int_0^x \sqrt{V(s)} \dd s\right|\,,\quad \varepsilon\in(0,1)\,,
	\end{equation}
	and gives
		\[\Re \mathscr{L}_\Phi - \lambda\geq \varepsilon V(x)-\lambda\,.\]
	This allows one to establish an estimate in the form
	\begin{equation}\label{eq.ellipticity0}
		\| e^{\frac \Phi h} u \| \leq C h^{- \kappa_1} \| e^{\frac \Phi h} ( \mathscr{L}-\lambda) u \| + C h^{-\kappa_2} \| u \|\,,
	\end{equation}
	for some powers $\kappa_1,\kappa_2>0$. The rather naive choice \eqref{eq.subAgmon} allows one to show in our specific situation that the groundstate eigenfunction is exponentially localized near $0$ at a scale of order $h^{\frac23}$ since, for $x$ small, we have $\Phi(x)\sim \frac{2\sqrt{1-\varepsilon}}{3}\sqrt{V'(0)}x^{\frac32}$.
	
	The choice \eqref{eq.subAgmon} was made in the case of non-degenerate wells, see \cite[Lemma 2.3]{HS} but it prevents one from taking the limiting weight $\Phi=\left|\int_0^x \sqrt{V(s)} \dd s\right|$, which is needed to get an exponentially good approximation of the eigenfunctions far from the well. To fix this issue  in \cite[Prop. 5.5]{HS}, the authors introduce a slightly better weight, which is a logarithmic correction of the previous one satisfying
	\[\Phi \simeq \left|\int_0^x \sqrt{V(s)}\mathrm{d}s\right| - Nh |\log h|\,,\] 
	with a large constant $N$. With their more involved choice, one can actually replace $\Phi$ by $\left|\int_0^x \sqrt{V(s)} \dd s\right|$ at the price of increasing the powers $\kappa_1$ and $\kappa_2$ by $N$. Note that smooth versions of such weights have recently been considered in a pseudo-differential framework in \cite[Section 4.2]{DR25}.
	
	In our situation with linear wells, none of these weights works since the eigenvalues are much larger (they are of order $\lambda \simeq h^{\frac23}$) than in the non-degenerate case (when they are of order $h$). For this reason, we rather consider the weight:
	\begin{equation}\label{eq.weight0}
		\Phi_h(x) = \sqrt{1- h} \int_0^x \sqrt{(V(s)- \lambda)_+} \dd s\,,
	\end{equation}
with which we can directly prove \eqref{eq.ellipticity0} with $\kappa_1 = \frac 43$, $\kappa_2= \frac 13$, and no fancy logarithmic correction is needed. Then, we can easily replace $\Phi_h$ by $S$ since $|\Phi_h - S| = \mathcal{O}(h)$ on any compact set.
	
	\begin{remark}
	 Our strategy also provides a simpler proof for non-degenerate wells, for which one can still use the weight \eqref{eq.weight0}, without any logarithmic correction. This choice allows one to prove simultaneously the optimal exponential decay of the eigenfunctions and their exponential approximation far from the well.
	\end{remark}

\subsection{Quasimodes and approximation} 
Thanks to the Agmon estimate \eqref{eq.ellipticity0} with the naive weight \eqref{eq.subAgmon}, we can see that most of the mass of the groundstate $f$ lives on an interval of size $h^{\frac 23}$ around the origin. Thus, a Taylor expansion near $x=0$ suggests considering the Airy operator  
\[
\mathscr{A} = -h^2 \partial_x^2 + V'(0)x\,.
\]
A first approximation of $f$ is then given by the groundstate of $\mathscr{A}$, \emph{i.e.}, uniformly on $[0,1]$,
\begin{equation}\label{eq.local.asymptotics}
	f(x) = c_0 h^{-\frac 13} {\rm Ai}( V'(0)^{\frac 13} h^{- \frac 23} x + z_1) + \mathcal{O}(h^{\frac 13}),
\end{equation}
where $\rm{Ai}$ denotes the Airy function, $z_1$ its first negative zero, and 
\[
c_0=\frac{V'(0)^{\frac16}}{\left(\int_0^{+\infty}\mathrm{Ai}^2(x+z_1)\mathrm{d}x\right)^{\frac12}} = \frac{V'(0)^{\frac{1}{6}}}{|\mathrm{Ai}'(z_1)|}
\]
is a normalizing constant. In $L^2$ norm, the error in \eqref{eq.local.asymptotics} is even $\mathcal{O}(h^{\tfrac{2}{3}})$. It also suggests that the first eigenvalues behave like in \eqref{eq.lambdajh}. These approximations are made rigorous in Proposition \ref{prop.L2approx}.

However \eqref{eq.local.asymptotics} is only good for $x$ of order $h^{\frac 23}$. Indeed, ${\rm{Ai}}$ decays exponentially fast at infinity, thus making the error term $\mathcal{O}(h^{\frac 13})$ too large to describe the announced behavior of $f$, when $x \gg h^{\frac 23}$, see Theorem \ref{thm.linear.wells}.

For larger values of $x$, one may want to consider a WKB Ansatz, \emph{i.e.}, a trial function of the form
\[ f^{\rm{WKB}} = e^{- \frac Sh} (a_0 + h a_1 + h^2 a_2 + \cdots ).\]
It is well-known that such an Ansatz can be built, and that it is a good approximation of the solution for large values of $x$ (i.e. far inside the classically forbidden region $\lbrace V(x) > \lambda \rbrace=\{x>x_\lambda\}$). However, this construction degenerates when we reach the classically allowed region, when $x$ is of order $h^{\frac 23}$. To overcome this problem, we give an alternative to the classical WKB method by changing the exponential profile of the WKB Ansatz, especially in the region $x \gtrsim h^{\frac 23}$, see Section \ref{sec.newquasimodes}. This Ansatz already appeared (for $N=0$) a long time ago in \cite[Section 5 (iii)]{O63}.

\begin{proposition}\label{thm.quasimode}
	Let $N \geq 0$, $K \geq 1$; let $\lambda=\lambda_1$. There exist $h_0,c >0$ and, for $h \in (0,h_0)$, a smooth function $\varphi$ on $[0,K]$ such that $\varphi'\geq c$ and  $\varphi'(x_\lambda)=V'(0)^{\frac13}+\mathcal{O}(h^{\tfrac{2}{3}})$ and such that the function
	\begin{equation*}
		f^{\rm{app}}_N(x) = \frac{\alpha_h}{\sqrt{\varphi'(x)}} {\rm Ai}(h^{-\frac 23} \varphi(x))\,,\quad\mbox{ with }\quad \alpha_h=f( x_\lambda)\frac{\sqrt{\varphi'( x_\lambda)}}{\mathrm{Ai}(h^{-\frac23}\varphi( x_\lambda))}\,,
	\end{equation*}
satisfies $f_N^{\mathrm{app}}(x_\lambda)=f(x_\lambda)$ and 
	\begin{equation}\label{eq.f0app}
		(\mathscr{L}-\lambda) f^{\rm app}_N = \mathcal{O}(h^{2N+2})f^{\rm app}_N,
		\end{equation}
	uniformly on $[x_\lambda,K]$. Moreover, 
	\[
	 \varphi (x) = \left( \frac 32 \int_{0}^x \sqrt{(V(s)-\lambda_1)_+} \dd s \right)^{\frac 23} + \mathcal{O}(h^2), 
	 \]
	for all $x\in [x_\lambda,K]$, with equality when $N=0$.
\end{proposition}

Compared to the WKB Ansatz, our construction is better suited for the shape of the well, through the use of the Airy function instead of the exponential. This idea could be generalized to other kinds of wells by using other special functions. \\
\begin{remark}
In fact, in order to prove Theorem \ref{thm.linear.wells}, we only need the construction with $N=0$. This is due to the low powers of $h^{-1}$ appearing in our Agmon estimates. For larger $N$, we get higher order approximations of $f$. See Section \ref{sec.proof}.
\end{remark}

\begin{remark}
\begin{enumerate}[--]
\item The function $\varphi$ that we construct satisfies $\varphi(0)=z_1 h^{\frac23}+\mathcal{O}(h^{\tfrac{4}{3}})$ and thus $f_N^{\mathrm{app}}$ almost satisfies the Dirichlet boundary condition. However, we will mainly use our Ansatz in the classically forbidden region.
\item In Lemma \ref{lem.almostC1}, it is shown that $(f^{\rm{app}}_N)'(x_\lambda)$ almost coincides with $f'(x_\lambda)$.
\item We recall that
\[\mathrm{Ai}(x)\underset{x\to +\infty}{=}\frac{1}{2x^{\frac14}\sqrt{\pi}}e^{-\frac23 x^{\frac32}}(1+o(1))\,.\]
Thus, for all fixed $x>0$, we have, when $h\to 0$,
\begin{equation}\label{eq.fapp.asymp}
	f^{\rm{app}}_N(x) = \frac{c}{h^{\frac16}V(x)^{\frac 14}} e^{- \frac 1h S(x)}\big( 1+ o(1) \big),
\end{equation}
with $c=\frac{c_0V'(0)^{\frac16}}{2\sqrt{\pi}}$ and where \eqref{eq.local.asymptotics} is used to estimate $\alpha_h$. This matches the asymptotic behavior in Theorem \ref{thm.linear.wells}.
\end{enumerate}
\end{remark}
In Section \ref{sec.def.psi}, the function $f_N^{\mathrm{app}}$ is used to construct an exponentially good quasimode defined as follows. We consider the function defined by
\[
\psi(x)=\begin{cases}
	f(x)&\mbox{ if }x\in[0, x_\lambda]\,,\\	
	\chi(x)\widetilde f_N^{\mathrm{app}}(x)&\mbox{ if } x_\lambda\leq x\,,
\end{cases}
\]
where $\chi$ is a smooth cutoff function and $\widetilde f_N^{\mathrm{app}}$ is a very small correction of $f^{\mathrm{app}}_N$ locally near $x_\lambda$ to make it $\mathscr{C}^1$ at $x_\lambda$ and belong to the domain of $\mathscr{L}$, see \eqref{eq.correction}. In Proposition \ref{prop.psi-f}, we prove that $\psi$ is a good $L^2$-approximation of $f$.

In Section \ref{sec.approximation}, we prove that $\psi$ is an exponentially good quasimode, see Lemma \ref{lem.expquasimode}. The crucial ingredient is \eqref{eq.f0app}, combined with the asymptotic behavior of the Airy function.

 Finally in Section \ref{sec.proof}, we use the Agmon estimates (apply Proposition \ref{prop.improvedAgmon} to $u=f-\psi$) and we get, for all fixed $x>0$,
\[e^{\Phi(x)/h}(f(x)-\psi(x))=\mathcal{O}(h^{2N})\,,\]
with a weight $\Phi$ such that $\Phi(x)=S(x)+\mathcal{O}(h)$.
Then, we have
\[f(x)=\psi(x)+\mathcal{O}(h^{2N})e^{-S(x)/h} \,.\]
For fixed $x$, and $h$ small enough, $\psi(x)=f_N^{\mathrm{app}}(x)$, which is estimated with \eqref{eq.fapp.asymp}, and thus Theorem \ref{thm.linear.wells} follows by taking $N=0$.

\section{Analysis of the single singular well}\label{sec.3}

	\subsection{Proof of the optimal Agmon estimates}\label{sec.improvedAgmon}

As explained in the introduction, the main idea to prove the Agmon estimates with optimal weight $S$ (as stated in Proposition \ref{prop.Ag}), is to first consider the quasi-optimal weight $\Phi$.

	\begin{proposition}\label{prop.improvedAgmon}
Let $E >0$. There exist $h_0$ and $C >0$ such that for all $h \in (0,h_0)$, all $0\leq \lambda \leq Eh^{\frac 23}$, the following holds. Let $M\geq 1$. Let us consider the bounded and Lipschitz function  
\[
\Phi(x)=\sqrt{1-h}\min\left(S(x),S\left(\frac{M}{2}\right)\right)\,,
\]  
with $S(x)=\int^x_{0}(V-\lambda)^{\frac12}_+$. For $u \in H^2(\R_+)\cap H^1_0(\R_+)$,
\begin{align*}
\| e^{\Phi / h} u \|	&\leq C h^{- \frac 43} \| e^{\Phi /h} (\mathscr{L}-\lambda)u \| + C h^{- \frac 13} \| u \|\,,\\
\|h\partial_x e^{\Phi / h} u \| 							&\leq C h^{- \frac 23} \| e^{\Phi /h} (\mathscr{L}-\lambda)u \| + C h^{ \frac 13} \| u \|\,.
\end{align*}
Additionally, there exists a constant $C_M>0$ so that 
\[
\|e^{\Phi / h}h^2\partial_x^2 u\|_{L^2(0,M/2)} 	\leq C_M h^{- \frac 43} \| e^{\Phi /h} (\mathscr{L}-\lambda)u \| + C_M h^{- \frac 13} \| u \|\,,
\]
\end{proposition}
	
\begin{proof}
By the Agmon formula (see \cite[Proposition 6.1]{DS}) or by using \eqref{eq.ReLPhi},
\begin{equation*}
	\Re \langle(\mathscr{L}-\lambda)u,e^{2\Phi/h}u\rangle=\|h\partial_x(e^{\Phi/h}u)\|^2+\int_{\R_+}\left(V-\lambda-\Phi'^2\right)|e^{\Phi/h}u|^2\mathrm{d}x\,.
\end{equation*}
The weight $\Phi$ is chosen so that $\Phi'^2 \leq (1-h)(V-\lambda)_+$. Thus,
\begin{equation*}
	\Re \langle(\mathscr{L}-\lambda)u,e^{2\Phi/h}u\rangle
	\geq\|h\partial_x(e^{\Phi/h}u)\|^2+\int_{\R_+}\left(h(V-\lambda)_+-(V-\lambda)_-\right)|e^{\Phi/h}u|^2\mathrm{d}x
\end{equation*}
Note that, since $\Phi = 0$ on $\lbrace V \leq \lambda \rbrace$ and $\lambda = \mathcal{O}(h^{\frac 23})$,
\begin{equation*}
	\int_{\R_+} (V-\lambda)_- | e^{\Phi/h} u |^2 \dd x \leq C h^{\frac 23} \|u \|^2.
\end{equation*}
Therefore,
\begin{equation*}
\| h \partial_x (e^{\Phi/h} u) \|^2 + h \int_{\R_+} (V-\lambda) | e^{\Phi/h} u |^2 \dd x \leq \Re \langle(\mathscr{L}-\lambda)u,e^{2\Phi/h}u\rangle + C h^{\frac 23} \| u \|^2.
\end{equation*}
The operator $- \frac {h^2} 2 \partial_x^2 + hV - h \lambda$ is bounded below by $c h^{\frac 43}$, with some positive constant $c$. Thus, after keeping half of the kinetic energy we find
\begin{equation}\label{eq.Agmon.qf}
\frac 12 \| h \partial_x (e^{\Phi/h} u) \|^2 + c h^{\frac 43} \| e^{\Phi/h} u \|^2 \leq \Re \langle(\mathscr{L}-\lambda)u,e^{2\Phi/h}u\rangle + C h^{\frac 23} \| u \|^2.
\end{equation}
Using the Young inequality we have
\begin{equation*}
{\rm Re} \langle ( \mathscr{L} - \lambda) u , e^{2\Phi/h} u \rangle \leq \frac{1}{2c} h^{-\frac 43} \|e^{\Phi/h}(\mathscr{L}-\lambda) u\|^2 + \frac{c}{2} h^{\frac 43} \| e^{\Phi /h} u \|^2,
\end{equation*}
and we deduce from \eqref{eq.Agmon.qf} the first two estimates:
\begin{equation*}
\|e^{\Phi/h}u\|\leq C h^{-\frac{4}{3}}\|e^{\Phi/h}(\mathscr{L}-\lambda) u\|+Ch^{-\frac13}\|u\|\,,
\end{equation*}
\begin{equation*}
\|h\partial_x(e^{\Phi/h}u)\|\leq C h^{-\frac 23}\|e^{\Phi/h}(\mathscr{L}-\lambda) u\|+Ch^{\frac13}\|u\|\,.
\end{equation*}
Finally, to get the second order estimates we use the first estimate and the triangle inequality,
\[
\|e^{\Phi/h}h^2\partial_x^2 u\|_{L^2(0,M/2)}\leq C\|e^{\Phi/h}(\mathscr{L}-\lambda)u\|+C_M\|e^{\Phi/h}u\|_{L^2(0,M/2)}.
\]
\end{proof}

To get the result as stated in Proposition \ref{prop.Ag}, we simply use Proposition \ref{prop.improvedAgmon} with $M = \max(1,2K)$, so that $S \leq (1+h) \Phi$ on $(0,K)$, and therefore
\[ \| e^{S / h} u \|_{L^2(0,K)} \leq C \| e^{\Phi / h} u \| \leq C h^{- \frac 43} \| e^{S /h} (\mathscr{L}-\lambda)u \| + C h^{- \frac 13} \| u \|,  \]
where in the last inequality we also used $\Phi \leq S$.
	
    \subsection{Linear approximation inside the well}
    
The Agmon estimates show that the groundstate $f$ of $\mathscr{L}$ is localized on an $h^{\frac 23}$ neighborhood of $0$. For this reason we can approximate $f$ by the groundstate $g_h^0$ of the Airy operator
\[
\mathscr{A} = -h^2 \partial_x^2 + V'(0)x\,.
\]
$g_h^0$ can be expressed in terms of the Airy function, solution to
\[
 - \partial_x^2 {\rm Ai} + x {\rm Ai} = 0, \quad \lim_{x \to + \infty} {\rm Ai}(x) = 0.
\]
The eigenvalue associated to $g_h^0$ is $|z_1| V'(0)^{\frac 23} h^{\frac 23}$, where $z_1$ is the first negative zero of ${\rm Ai}$. Higher order Taylor expansions of $V$ yield good approximations of $f$ in $L^2$ and even $\mathscr{C}^1$ norm.
\begin{proposition}\label{prop.L2approx}
For all $k\geq 1$, we have 
\[
\lambda_k(h)=|z_k|V'(0)^{\frac23}h^{\frac23}+\mathcal{O}(h^{\frac43})\,,
\]
where $z_k$ is the $k$-th negative zero of the Airy function. Moreover, there exists a family of functions in the Schwartz class $(g_j)_{j\geq 0}$, all vanishing at $0$, the first one being given by
\[
g_{0}(x)=c_0\mathrm{Ai}(V'(0)^{\frac13}x+z_1)\,,
\]
with $c_0=V'(0)^{\frac16}|\mathrm{Ai}'(z_1)|^{-1}$, such that, for $N \geq 0$, the function 
\[
g^N_h(x)=h^{-\frac13}\sum_{j=0}^N h^{\frac{2j}{3}}g_j(h^{-\frac23}x)\,,
\]
satisfies
\[
\|f-\chi g_h^N\|_{L^2(\R_+)}=\mathcal{O}(h^{\frac{2N+2}{3}}) \quad {\rm and} \quad \|f-\chi g_h^N\|_{\mathscr{C}^1(\R_+)}=\mathcal{O}(h^{\frac{2N}{3}-\frac{1}{2}})\,. 
\]
for any cutoff function $\chi$ supported and equal to $1$ near $0$. For some such cutoff, set $g_h = \chi g_h^N$. 
\end{proposition}
    
    \begin{proof}
The proof follows classical arguments. Let us recall the main steps. We start by constructing a rather naive quasimode. By rescaling, the operator $\mathscr{L}$ is unitarily equivalent to $\mathscr{L}^{\mathrm{resc}}=-h^{\frac23}\partial_y^2+V(h^{\frac23}y)$ which we can formally Taylor expand as		
\[
\mathscr{L}^{\mathrm{resc}}\sim h^{\frac23}\sum_{j\geq 0}h^{\frac{2j}{3}}\mathscr{L}_j\,,\quad \mathscr{L}_0 = - \partial_y^2 + V'(0) y\,,\quad \mathscr{L}_j=\frac{V^{(j+1)}(0)}{j+1!}y^{j+1}\,(j\geq 1),
\]
This suggests looking for a quasi-eigenvalue $\mu=h^{\frac23}\sum_{j=0}^N h^{\frac{2j}{3}}\mu_j$ and for a quasimode in the form
\[
g(y)=\sum_{j=0}^N h^{\frac{2j}{3}}g_j(y)\,,
\]
Solving formally 		
\[
(\mathscr{L}^{\mathrm{resc}}-\mu)g=0
\]
gives a sequence of equations. The first one is
\[
(\mathscr{L}_0-\mu_0)g_0=0\,,
\]
which leads to the choice $\mu_0=|z_1|V'(0)^{\frac23}$ (the groundstate eigenvalue of the Airy operator $\mathscr{L}_0$) and to
\[
g_{0}(y)=c_0\mathrm{Ai}(V'(0)^{\frac13}y+z_1)\,,
\]
the corresponding normalized eigenfunction. 
    	
Then, the next $(\mu_j,g_j)$ are determined by an induction procedure requiring us to solve an equation in the form
\begin{equation}\label{eq.gj}
(\mathscr{L}_0-\mu_0)g_j=\mu_j g_0+R_j\,
\end{equation}
where $R_j$ is known and belongs to the Schwartz class. Since the left-hand side must be orthogonal to $g_0$, we need to choose $\mu_j$ such that the right-hand side also is, i.e. $\mu_j = - \langle R_j, g_0 \rangle$. Then, for such right-hand side, there is a unique solution $g_j$ to \eqref{eq.gj} so that 
\[
\Big\|\sum_{\ell=0}^j h^{\frac{2\ell}{3}} g_\ell \Big\|^2 = 1 + \mathcal{O}(h^{\frac{2(j+1)}{3}}). 
\]
This determines the sequence $\mu_j$ and $g_j$. The $g_j$ are Schwartz class, this is shown by variation of parameters. For any $C>0$, using the Taylor formula with remainder to order $N+ 2$ for $V$, and the fact that the $g_j$'s are Schwartz class, we get
\[
(\mathscr{L}^{\mathrm{resc}} - \mu)\sum_{j=0}^N h^{\frac{2j}{3}} g_j = \mathcal{O}_{L^2([0,C h^{-2/3}])}(h^\frac{2N+4}{3}). 
\]
(We have to restrict to $C<\infty$ because we have not made any assumption on the size of the derivatives of $V$ at infinity). 
    	
Rescaling back, and using the fact $[\mathscr{L},\chi]$ is supported where $g_h^N$ is $\mathcal{O}(h^\infty)$, we deduce that
\begin{equation}\label{eq.ghquasimode}
\|(\mathscr{L}-\mu)(\chi g_h^N)\|_{L^2(\R_+)}=\mathcal{O}(h^{\frac{2N+4}{3}})\,.
\end{equation}
Thanks to the spectral theorem, this shows that $\mathrm{dist}(\mathrm{sp}(\mathscr{L}),\mu)\leq Ch^{\frac{2N+4}{3}}$. Playing the same game by choosing other eigenfunctions of the Airy operator associated with the next zeroes, we build a sequence of quasi-eigenvalues, all of order $\mathcal{O}(h^{\frac23})$.
    	
By means of Proposition \ref{prop.improvedAgmon}, it can be established that there exists $c>0$ such that, for all $F\in\mathscr{E}:=\ker(\mathscr{L}-\lambda_1(h))\oplus\ker(\mathscr{L}-\lambda_2(h))$,
\begin{equation}\label{eq.Agmonh23}
\int_{\R_+}e^{\frac{c x}{h^{2/3}}}|F(x)|^2\mathrm{d}x\leq C\|F\|^2\,.
\end{equation}
To prove this, one can proceed directly with the usual strategy for Agmon estimate. Considering $\mathscr{Q}$ the quadratic form associated with $\mathscr{L}$, we have
\[
\mathscr{Q}(F)= \int_{\R_+}(h^2|F'(x)|^2+V(x)|F(x)|^2)\mathrm{d}x\leq \lambda_2(h)\|F\|^2
\]
and then, by using the Taylor expansion of $V$ and estimating the remainders thanks to \eqref{eq.Agmonh23},
\[
\mathscr{Q}(F)\geq  \int_{\R_+}(h^2|F'(x)|^2+V'(0)x|F(x)|^2)\mathrm{d}x-Ch^{\frac43}\|F\|^2\,.
\]
By the min-max theorem, this shows that
\[
\lambda_2(h)\geq |z_2|V'(0)^{\frac23}h^{\frac23}-Ch^{\frac43}\,,
\]
and therefore $|\lambda_1(h)-\mu|={\rm{dist}}({\rm sp}( \mathscr{L}),\mu)=\mathcal{O}(h^{\frac{2N+4}{3}})$. (The same argument gives similar estimation for the next eigenvalues that are thus all simple). With the spectral theorem and \eqref{eq.ghquasimode}, we deduce that
\[
\|f-g_h\|_{L^2(\R_+)}=\mathcal{O}(h^{\frac{2N+2}{3}})\,.
\]
To obtain a uniform estimate, we can use that for $u\in H^2\cap H^1_0$ supported in $[0,K]$,
\begin{align*}
\|(\mathscr{L}-\mu)u\|^2 	&= \| h^2 u''\|^2 + \|(V-\mu)u\|^2 + 2h^2 \Re \langle (V-\mu) u', u'\rangle - h^2\langle V'' u,u\rangle. \\
h^2\|u'\|^2 				&\leq \langle(\mathscr{L}-\mu)u,u\rangle + \mu\|u\|^2.
\end{align*}
It follows that 
\[
\| h^2 \partial_x^2 u \|_{L^2}^2 \leq C h^{-2/3} \| (\mathscr{L} - \mu)u \|_{L^2}^2 + C h^{2/3} \|u  \|_{L^2}^2\,,
\]
with constants depending on $K$ (actually $\|V''\|_{L^\infty(0,K)}$). Inserting a cutoff against $f-g_h$, and using that $f$ is $\mathcal{O}(h^\infty)$ for $x>1$, we deduce that 
\[
\| h^2 \partial_x^2(f-g_h)\|_{L^2(]0,K[)} \leq C h^{\frac{2N+3}{3}}.
\]
For large values of $x>K$, we use the fact that $V> \lambda$ and $f>0$, so that 
\[
f'' = h^{-2}(V-\lambda)f >0. 
\]
This implies that $f'$ is non-decreasing, but it goes to $0$, so $f'(K) \leq f'(x) \leq 0$, and $0\leq f(x) \leq f(K)$. 
\end{proof}
    
    Since ${\rm Ai}$ is exponentially decaying at infinity, the approximation of $f$ given by Proposition \ref{prop.L2approx} is only good on an $h^{2/3}$ neighborhood of $0$. Still, this allows to compute the asymptotic value of $f$ at the boundary of the well, i.e. the point $x_\lambda$ such that $V(x_\lambda) = \lambda_1$.
    
\begin{corollary}\label{cor.fxlambda}
We have
\[
f(x_\lambda)\underset{h\to 0}{\sim} c_0h^{-\frac13}\mathrm{Ai}(0)\,.
\]
\end{corollary}
    
\begin{proof}
First note that $x_\lambda$ must be of order $h^{\frac 23}$. By a Taylor approximation we deduce that
\[
\lambda_1 = V(x_\lambda) = V'(0) x_\lambda + \mathcal{O}(h^{\frac 43}) ,
\]
and since $\lambda_1 = |z_1| V'(0)^{\frac 23} h^{\frac 23} + \mathcal{O}(h^{\frac 43})$ we obtain
\begin{equation}\label{eq.xlambda} 
x_\lambda = |z_1| V'(0)^{-\frac 13} h^{\frac 23} + \mathcal{O}(h^{\frac 43}). 
\end{equation}
Using the uniform estimates in Proposition \ref{prop.L2approx}, we find that
\begin{align*}
f(x_\lambda) 	&= h^{-\frac 13} g_0(h^{-\frac 23} x_\lambda) + \mathcal{O}(h^{\frac 13})\\
				&= c_0 h^{-\frac 13} \mathrm{Ai}(V'(0)^{\frac 13} h^{-\frac 23} x_\lambda + z_1) + \mathcal{O}(h^{\frac 13}),
\end{align*}
and the result follows using \eqref{eq.xlambda}.
\end{proof}

\subsection{A construction of quasimodes}\label{sec.newquasimodes}
    	We now construct a quasimode that will be exponentially good. It will be shown later to be a correct approximation of the true groundstate $f$ at large and intermediate distance from the well. 
    	
    	In the following, we denote by $\lambda = \lambda_1$ the first eigenvalue of $\mathscr{L}$. We recall that $x_\lambda$ denotes the boundary of the well, i.e. the solution to $V(x_\lambda) = \lambda$.  It will be convenient to extend $V$ to a smooth function on $(-\varepsilon, \infty)$, for some small $\varepsilon >0$. 
    	
We introduce the weight
\begin{equation}\label{def.phi0}
\varphi_0(x)=\begin{cases}
			\left(\frac32\int_{x_\lambda}^x\sqrt{V-\lambda}\right)^{\frac23}&\mbox{ if } x\geq x_\lambda\,,\\
			-\left(\frac32\int_{x}^{x_\lambda}\sqrt{\lambda-V}\right)^{\frac23}&\mbox{ if } x \leq x_\lambda\,.
\end{cases}
\end{equation}
In fact, $\varphi_0$ defines a smooth function.
	
\begin{lemma}\label{lem.phi.regularity}
Let $M\geq 1$. There exists $\varepsilon>0$ such that the function $\varphi_0$ is smooth and increasing on $(-\varepsilon,M)$. Moreover, there exist $c>0$ and $h_0>0$ such that, for all $h\in(0,h_0)$ and all $x\in(-\varepsilon,M)$, we have $\varphi_0'(x)\geq c>0$. In addition, we have
\begin{equation}\label{eq.phi0prime}
V(x)-\lambda=\varphi_0(x)\varphi'_0(x)^2\,.
\end{equation}	
\end{lemma}
	
\begin{proof}
The regularity problem appears near $x=0$. Thanks to the Taylor formula, we get
\[
V(x)-V(x_\lambda)=(x-x_\lambda)\int_0^1 V'(x_\lambda+t(x-x_\lambda))\mathrm{d}t\,.
\]
For $h$ small enough and $x$ in a neighborhood of $0$, we can write
\[
\int_0^1 V'(x_\lambda+t(x-x_\lambda))\mathrm{d}t=R(x-x_\lambda)^2\,,
\]
for some smooth function $R$ such that $R\geq c_0>0$. Then, for $x\geq x_\lambda$, we have
\[
\int_{x_\lambda}^x\sqrt{V(s)-\lambda}\,\mathrm{d}s=\int_0^{x-x_\lambda}\sqrt{\theta}R(\theta)\,\mathrm{d}\theta=\frac23(x-x_\lambda)^{\frac32}\int_0^1 R(v^{\frac23}(x-x_\lambda))\,\mathrm{d}v
\]
so that
\[
\varphi_0(x)=(x-x_\lambda)\left(\int_0^1 R(v^{\frac23}(x-x_\lambda))\,\mathrm{d}v\right)^{\frac23}\,.
\]
We can check that this formula also holds for $x\leq x_\lambda$. We see that $\varphi_0$ is a smooth function of $x$, since $R \geq c_0 >0$ everywhere. Finally we can compute $\varphi_0'$ using \eqref{def.phi0} and we immediately see that $\varphi_0$ is increasing and satisfies \eqref{eq.phi0prime}.
\end{proof}	
	
The next proposition gives the construction of our quasimode, which is a combination of the Airy function (well-suited for the shape of the well) and an iterative perturbation $\varphi$ of $\varphi_0$.
	
\begin{proposition}\label{prop.quasimode0}
Let $K\geq 1$. There exist $\varepsilon_0, h_0>0$ and, for $h \in (0,h_0)$, smooth increasing functions $\varphi_1,\ldots,\varphi_N,\dots$ defined on $\varphi_0((-2\varepsilon_0,K+\varepsilon_0))$ such that the following holds. Uniformly and for all $j\geq 1$ we have
\[
\varphi_j(x)=x+\mathcal{O}_{C^\infty}(h^{2j})\,.
\]
Thus, for $N>0$ and $h>0$ small enough, $\varphi=\varphi_N\circ\ldots\circ\varphi_1\circ\varphi_0$ is well defined and we let
\[
f_N(x)=\frac{1}{\sqrt{\varphi'(x)}}\mathrm{Ai}(h^{-\frac23}\varphi(x))\,.
\]
Then, uniformly on $(-\varepsilon_0,K)$ we have
\[
(\mathscr{L}-\lambda)f_N=\mathcal{O}(h^{2N+2})f_N\,.
\]
\end{proposition}
	
\begin{proof}
To construct the announced quasimode, we first look for it in the form 
\[
u(x)=\alpha(x)v(\varphi(x))\,.
\]
The second derivative of such a function is given by
\[
u''(x)=\alpha'' v+2\alpha'\varphi'v'+\alpha\varphi'' v'+\alpha\varphi'^2v''
\]
and we choose
\[
\alpha(x)=\frac{1}{\sqrt{\varphi'(x)}}\,,
\]
in order to kill first order derivatives in $v$. Letting $y=\varphi(x)$, a computation gives
\[
(\mathscr{L}-\lambda)u=\alpha(x)\varphi'(x)^2\left(-h^2  v''(y)+\frac{V(\varphi^{-1}(y))-\lambda}{\varphi'(\varphi^{-1}(y))^2}v(y)-h^2\frac{\alpha''}{\alpha[\varphi'(\varphi^{-1}(y))]^2} v(y)\right)\,.
\]
By choosing $\varphi(x)=\varphi_0(x)$ for all $x\in(-\varepsilon,M)$, we get
\[
\frac{V(\varphi^{-1}(y))-\lambda}{\varphi'(\varphi^{-1}(y))^2}=y\,,
\]
for all $y\in(\varphi_0(-\varepsilon),\varphi_0(M))$. We let the remainder term be
\[
r_0=\Big[ (\varphi_0')^{-3/2} ((\varphi_0')^{-1/2})''\Big]\circ\varphi_0^{-1}\,,
\]
so that
\[
 (\mathscr{L} - \lambda)u = (\varphi_0')^{\frac 32} \left( -h^2 \partial_y^2 v + y v - h^2 r_0 v \right).
 \]
If $r_0$ were $0$ we could choose $v={\rm Ai}(h^{-2/3}y)$. Instead we will recursively reduce the size of the remainder term. We look for $v$ in the form
\[
v(y)=\frac{1}{\sqrt{\varphi'_1(y)}}w(\varphi_1(y))\,.
\]
The map $y\mapsto y-h^2r_0(y)=W_1(y)$ is increasing when $h$ is small enough. There exists a unique $y_\lambda$ such that $y_\lambda-h^2r_0(y_\lambda)=0$. We let, for all $y\in(\varphi_0(-\varepsilon),\varphi_0(M))$,
\[
\varphi_1(y)=\begin{cases}
	\left(\frac32\int_{y_\lambda}^y\sqrt{W_1(s)}\mathrm{d}s\right)^{\frac23}&\mbox{ if } y\geq y_\lambda\,,\\
	-	\left(\frac32\int_{y}^{y_\lambda}\sqrt{-W_1(s)}\mathrm{d}s\right)^{\frac23}&\mbox{ if } y\leq y_\lambda\,.
\end{cases}
\]
The function $\varphi_1$ is smooth, increasing, and we have, uniformly in $y\in(\varphi_0(-\varepsilon),\varphi_0(M))$,
\[
\varphi'_1(y)=1+\mathcal{O}(h^2)\,.
\]

In the same way as previously, we get
\[
(-h^2\partial^2_y+y-h^2r_0(y))v=\varphi'_1(y)^{\frac32}\left(-h^2\partial^2_z+z-h^2r_1(z)\right)w\,,
\]
with $z = \varphi_1(y)$ and
\[
r_1=\Big[(\varphi_1')^{-3/2}((\varphi_1')^{-1/2})''\Big]\circ\varphi_1^{-1}=\mathcal{O}(h^2)\,.
\]	
This procedure can be continued by induction: at step $N$ we obtain $\varphi = \varphi_N\circ\ldots\circ\varphi_1\circ\varphi_0$, and 
\[
 (\mathscr{L} - \lambda)u = (\varphi')^{\frac 32} \left( -h^2 \partial_z^2  + z + \mathcal{O}(h^{2N+2}) \right) w \circ \varphi. 
 \]
We end the induction by taking $w= {\rm{Ai}}(h^{-2/3}z)$, which gives 
\[
 (\mathscr{L} - \lambda)u = \mathcal{O}(h^{2N+2}) u.
\]
	\end{proof}

\subsection{Approximation of the groundstate in $L^2$ norm}\label{sec.def.psi}

Proposition \ref{prop.quasimode0} suggests considering the function given, for all $K>x\geq x_\lambda$, by 
\begin{equation}\label{eq.fNapp}
f_N^{\mathrm{app}}(x)=\frac{\alpha_h}{\sqrt{\varphi'(x)}}\mathrm{Ai}(h^{-\frac23}\varphi(x))\,,
\end{equation}
where $\varphi=\varphi_N\circ\ldots\circ\varphi_1\circ\varphi_0$, and with normalization constant 
\[
\alpha_h=f( x_\lambda)\frac{\sqrt{\varphi'( x_\lambda)}}{\mathrm{Ai}(h^{-\frac23}\varphi( x_\lambda))}.
\]
The constant $\alpha_h$ is chosen so that $f^{\mathrm{app}}_N(x_\lambda)=f(x_\lambda)$, and we find that
\[
\alpha_h \sim c_0 h^{-\frac{1}{3}} V'(0)^{\frac{1}{6}}. 
\]
However, the derivatives do not match at $x=x_\lambda$. In fact, $\left(f_N^{\mathrm{app}}\right)'(x_\lambda)$ almost equals $f'(x_\lambda)$.

\begin{lemma}\label{lem.almostC1}
We have
\[
(f-f^{\mathrm{app}}_N)'( x_\lambda)=\mathcal{O}(h^{2N+ \frac 13})\,.
\]
\end{lemma}
\begin{proof}
We consider the Wronskian
\[
\mathcal{W}=f(f_N^{\mathrm{app}})'-f'f_N^{\mathrm{app}}
\]
which satisfies
\[
-h^2\mathcal{W}'(x)=f(x)(\mathscr{L}-\lambda)f_N^{\mathrm{app}}\,.
\]
We get
\[
f( x_\lambda)(f_N^{\mathrm{app}}-f)'( x_\lambda)=\mathcal{W}( x_\lambda)=\mathcal{W}(K)+\int_{ x_\lambda}^K h^{-2}f(\mathscr{L}-\lambda)  f_N^{\mathrm{app}}\,.
\]
The conclusion follows by using that $\mathcal{W}(K)=\mathcal{O}(h^{\infty})$ and by recalling Corollary \ref{cor.fxlambda}.
\end{proof}

We make a small perturbation of $f_N^{\rm app}$ in order to correct its derivative at $x=x_\lambda$. In what follows, we consider a smooth non increasing function $\chi$ such that $0\leq\chi\leq 1$, equal to $1$ on $[0,\frac{3M}{4})$ and $0$ on $(M,+\infty)$. Let us add a small perturbation of $f_N^{\mathrm{app}}$ by considering
\begin{equation}\label{eq.correction}
\widetilde f_N^{\mathrm{app}}=f_N^{\mathrm{app}}+(x-x_\lambda)\beta_h\chi_h(x)\,,
\end{equation}
where $\beta_h=(f-f_N^{\mathrm{app}})'( x_\lambda)$ and $\chi_h(x)=\chi(h^{-\frac23}x)$. Provided $M$ is large enough, $\chi_h(x_\lambda)=1$, and this function has derivative $(\widetilde f_N^{\rm app})'(x_\lambda) = f'(x_\lambda)$ at the boundary of the well.

\begin{proposition}\label{prop.psi-f}
We let
\[
\psi(x)=\begin{cases}
	f(x)&\mbox{ if }x\in[0, x_\lambda]\,,\\	
	\chi(x)\widetilde f_N^{\mathrm{app}}(x)&\mbox{ if } x\geq x_\lambda\,,
\end{cases}
\]
which is of class $\mathscr{C}^1$ and belongs to the domain of $\mathscr{L}$ by construction. It satisfies in $L^2$
\[
(\mathscr{L}-\lambda)\psi=\mathcal{O}(h^{2N+\frac 53})\|\psi\|\,.
\]
Moreover,
\[
\|\psi- f\|_{L^2(\R_+)}=\mathcal{O}(h^{2N+1})\,.
\]
\end{proposition}

\begin{proof}
We estimate the two terms separately in \eqref{eq.correction} to find
\[
\|(\mathscr{L}-\lambda)\psi\| = \mathcal{O}(h^{2N+2} \alpha_h) + \mathcal{O}(h^{\frac 43} \beta_h) = \mathcal{O}(h^{2N+\frac 53}).
\]
Using the same arguments as in the proof of Proposition \ref{prop.L2approx}, the spectral theorem gives that
\[
\|\psi-\|\psi\| f\|_{L^2(\R_+)}=\mathcal{O}(h^{2N+1})\,.
\]
Then, we notice that
\[
|1-\|\psi\||\|f\|_{L^2(0,x_\lambda)}=\|f-\|\psi\| f\|_{L^2(0,x_\lambda)}\leq \|\psi-\|\psi\| f\|_{L^2(\R_+)}\,.
\]
We use Proposition \ref{prop.L2approx} to see that $\|f\|_{L^2(0,x_\lambda)}\geq c>0$. Thus,
\[
|1-\|\psi\||=\mathcal{O}(h^{2N+1})\,.
\]
\end{proof}

\subsection{Approximation of the groundstate in exponentially weighted norm}\label{sec.approximation}

In this section we improve the estimates on our quasimode $\psi$, defined in Proposition \ref{prop.psi-f}. We show that $\psi$ is a good approximation of the true groundstate $f$ in exponentially weighted norm. 

\begin{lemma}\label{lem.expquasimode}
We have	
\[
e^{\Phi/h}(\mathscr{L}-\lambda)\psi=\mathcal{O}(h^{2N+\frac53})\,,
\]
where $\Phi$ is defined in Proposition \ref{prop.improvedAgmon}.
\end{lemma}

\begin{proof}
We only need to consider the region $x \geq x_\lambda$, where $\psi = \chi \widetilde f_N^{\rm app}$. Then,
\[
\begin{split}
e^{\Phi/h}(\mathscr{L}-\lambda)( \chi \widetilde f_N^{\mathrm{app}})	&=e^{\Phi/h}[-h^2\partial_x^2,\chi]\widetilde f_N^{\mathrm{app}}+e^{\Phi/h}\chi(\mathscr{L}-\lambda)\widetilde f_N^{\mathrm{app}}\\
																		&=e^{\Phi/h}[-h^2\partial_x^2,\chi] f_N^{\mathrm{app}}+e^{\Phi/h}\chi(\mathscr{L}-\lambda)(f^{\mathrm{app}}_N+\beta_h(x-x_\lambda)\chi_h)\,,
\end{split}
\]
where the first `tilde' was removed by support considerations.
The last term can be estimated using that $\Phi \leq Ch$ on the support of $\chi_h$.
\[
e^{\Phi/h}(\mathscr{L}-\lambda)(\chi \widetilde f_N^{\mathrm{app}})=e^{\Phi/h}[-h^2\partial_x^2,\chi] f_N^{\mathrm{app}}+e^{\Phi/h}\chi(\mathscr{L}-\lambda) f_N^{\mathrm{app}}+\mathcal{O}(h^{\frac 43}|\beta_h|)\,.
\]
Let us notice that
\[
\|e^{\Phi/h}\chi(\mathscr{L}-\lambda) f_N^{\mathrm{app}}\|\leq Ch^{2N+2}\|\chi e^{\Phi/h}f_N^{\mathrm{app}}\|\leq Ch^{2N+\frac53}\|\chi e^{\Phi/h}\mathrm{Ai}(h^{-\frac23}\varphi(x))\| \,,
\]
where we used Proposition \ref{prop.quasimode0} and Corollary \ref{cor.fxlambda}. Then bound the Airy function, for $z\in\R$,
\[
|\mathrm{Ai}(z)|\leq Ce^{-\frac{2}{3}z_+^{\frac32}},
\]
to get that
\[
\|\chi e^{\Phi/h}\mathrm{Ai}(h^{-\frac23}\varphi(x))\|\leq C\|\chi e^{\Phi/h}e^{-\frac23\varphi_+^{\frac32}/h}\|\leq C\|\chi e^{S/h}e^{-\frac23\varphi_+^{\frac32}/h}\|
\]
where we used that $\Phi\leq S$. Then, by construction we have, for all $K>x\geq x_\lambda$,
\[
\varphi(x)=\left(\frac{3}{2}\int_{x_\lambda}^x\sqrt{V-\lambda}\right)^{\frac23}+\mathcal{O}(h^2)\,,
\]
and by splitting the integral into two parts $x\geq x_\lambda$ and $x\leq x_\lambda$, we deduce
\[
\|\chi e^{\Phi/h}\mathrm{Ai}(h^{-\frac23}\varphi(x))\|\leq C\,.
\]
There remains to estimate the commutator 
\[
e^{\Phi/h}[-h^2\partial_x^2,\chi]f^{\mathrm{app}}_N=-h^2\chi'' e^{\Phi/h}f_N^{\mathrm{app}}-2h^2\chi' e^{\Phi/h}(f_N^{\mathrm{app}})'\,,
\]
which is supported in $(\frac{3M}{4},M)$.
The first term can be estimated by 
\[
\begin{split}
\|-h^2\chi'' e^{\Phi/h}f_N^{\mathrm{app}}\|	&\leq Ch^2 \|e^{\Phi/h}\mathrm{Ai}(h^{-\frac23}\varphi)\|_{L^2(\frac{3M}{4},M)}\\
											&\leq Ch^2 \|e^{\Phi/h}e^{-\frac{2}{3h}\varphi^{\frac32}}\|_{L^2(\frac{3M}{4},M)}\\
											&\leq Ch^2\|e^{(S(M/2)-S)/h}\|_{L^2(\frac{3M}{4},M)}\,.	 
\end{split}
\]
Thus,
\[
\|-h^2\chi'' e^{\Phi/h}f_N^{\mathrm{app}}\|\leq Ch^2e^{(S(M/2)-S(3M/4))/h}=\mathcal{O}(h^{\infty})\,,
\]
and we can deal with the other term in a similar way. We deduce that
\[
e^{\Phi/h}(\mathscr{L}-\lambda)(\chi \widetilde f^{\mathrm{app}}_N)=\mathcal{O}(h^{2N+\frac53})\,.
\]
Since $\psi$ is $\mathscr{C}^1$ and $(\mathscr{L}-\lambda)f=0$, the result follows.
\end{proof}

\subsection{Proof of Theorem \ref{thm.linear.wells}}\label{sec.proof}
We apply the Agmon estimates from Proposition \ref{prop.improvedAgmon} to the function $u=\psi- f$. We use Lemma \ref{lem.expquasimode} and Proposition \ref{prop.psi-f} to estimate the right-hand-side. 
\begin{equation}\label{eq.estimate.H2}
\| e^{\Phi / h} u \|+\|(e^{\Phi / h}h^2\partial_x^2 u)\|_{L^2(0,M/2)} = \mathcal{O}(h^{2N+ \frac 13}).
\end{equation}
and
\begin{equation}\label{eq.estimate.H1}
\|h\partial_x(e^{\Phi / h} u)\|= \mathcal{O}(h^{2N+ 1}).
\end{equation}
Then, by the Agmon inequality we can estimate the supremum norm,
\[
\| e^{\Phi /h} (\psi-f) \|_\infty  \leq C \| e^{\Phi /h} (\psi-f) \|_{L^2}^{\frac 12} \|  e^{\Phi /h} (\psi- f) \|_{H^1}^{\frac 12} = \mathcal{O}(h^{2N}) .
\]
Therefore, given a fixed $x>0$ we have $f(x)=f_N^{\mathrm{app}}(x)+\mathcal{O}(h^{2N})e^{-\Phi(x)/h}$, with 
\[
f_N^{\mathrm{app}}(x)=\frac{\alpha_h}{\sqrt{\varphi'(x)}}\mathrm{Ai}(h^{-\frac23}\varphi(x))\sim \frac{c_0 V'(0)^{\frac16}}{2\sqrt{\pi}}\frac{h^{-\frac16}}{\sqrt{\varphi'(x)}\varphi(x)^{\frac14}}e^{-\frac{2}{3h}\varphi(x)^{\frac32}}\,,
\]
so that
\[
f_N^{\mathrm{app}}(x)\sim \frac{c_0 V'(0)^{\frac16}}{2\sqrt{\pi}}\frac{h^{-\frac16}}{\sqrt{\varphi'_0(x)}\varphi_0(x)^{\frac14}}e^{-\frac{2}{3h}\varphi_0(x)^{\frac32}}=\frac{c_0V'(0)^{\frac16}}{2\sqrt{\pi}}\frac{h^{-\frac16}}{\sqrt{\varphi'_0(x)}\varphi_0(x)^{\frac14}}e^{-S(x)/h}\,.
\]
This  finishes the proof of Theorem \ref{thm.linear.wells} since
\[
\varphi_0(x)^{\frac14}=\left(\frac32\int^x_{x_\lambda} \sqrt{V-\lambda}\right)^{\frac16}\,,\quad
\sqrt{\varphi_0'(x)}=\left(\frac 32\int^x_{x_\lambda} \sqrt{V-\lambda}\right)^{-\frac16}(V(x)-\lambda)^{\frac14}\,.
\]
Note that $N=0$ is enough in the quasimode construction, for the error terms to be smaller than the leading order.

\begin{remark}
Note that we can similarly use \eqref{eq.estimate.H2} to estimate the derivative of $f$, on the interval $(0,M/2)$:
\begin{align*}
 \| e^{\Phi /h} (\psi'-f') \|_\infty \leq C \| e^{\Phi /h} (\psi'-f') \|_{L^2}^{\frac 12} \|  e^{\Phi /h} (\psi'- f') \|_{H^1}^{\frac 12} = \mathcal{O}(h^{2N-1}),
\end{align*}
using the control on the $H^2$ norm.
\end{remark}

\section{Applications to tunneling}\label{sec.4}
	
\subsection{Double linear wells}

Once we have proven sharp exponential decay estimates on the single-well eigenfunctions, together with pointwise asymptotic formulas, we have all the necessary ingredients to obtain tunneling formulas for the double-well problem. 

In this section we briefly explain how to obtain a tunneling formula between linear wells using our single well results. We will only recall the main ideas, since they follow a standard strategy.

Let us consider a smooth non-negative electric potential $W$ which is an even function on $I=[-a,a]$. We assume that $W$ has exactly two minima at $\pm a$ such that $W(\pm a)=0$ and $W'(a)=-W'(-a)<0$. We consider the operator
\[
\mathscr{L}=-h^2\partial_x^2+W\,,
\]
acting on the domain $H^2(I)\cap H^1_0(I)$. A tunneling effect between the two wells yields the following spectral gap.

\begin{proposition}\label{prop.application1}
There exists a universal constant $\kappa>0$ such that
\[
\lambda_2(h)-\lambda_1(h)=\kappa h^{\frac23}|W'(a)|^{\frac23}e^{-S_h/h}(1+o(1))\,,\quad S_h=\int_{-a}^a\sqrt{\big(W-h^{\frac23}|z_1|| W'(a)|^{\frac23}\big)_+}\,.
\]	
\end{proposition}
The tunneling estimate of Proposition \ref{prop.application1} does not seem to be already known. To the authors' knowledge, the only similar result appears in \cite[Section 1.5]{popoff} in the special case when $W(x)=1-|x|$ and $a=1$. There the author describes the eigenfunctions explicitly by means of the Airy functions $\mathrm{Ai}$ and $\mathrm{Bi}$. Note that, in this case,
\[
S_h=\frac{4}{3}-2|z_1|h^{\frac23}+\mathcal{O}(h)\,.
\]
Of course, we cannot directly apply our result due to the peak at $x=0$, but our estimate agrees with \cite[Corollaire 1.18]{popoff}. Note also that a similar formula is obtained in \cite[Corollaire 3.2]{MR88} in the case of degenerate wells (with a behavior $x^{2\sigma}$ near the minima, with $\sigma\geq 1$)\footnote{The correct prefactor in \cite[Corollaire 3.2]{MR88} involves $\mu^{1/(\sigma+1)}$ and not $\mu^{2/(\sigma+1)}$ which comes from a mistake on page 1182.} and it is consistent with ours if one were allowed to take $\sigma=\frac12$, but their Agmon distance does not include a perturbation like our $h^{\frac23}$.

\begin{proof}
Let us only sketch the proof by explaining how to use our estimate of the one-well groundstate in Theorem \ref{thm.linear.wells}. The details of the computations can be found in \cite[Chapter 6]{DS} and also in recent articles, see for instance \cite[Section 5]{DR25} or \cite[Section 5]{AFHR25}.

Given $a/2>\eta>0$, we consider a smooth single-well potential $W_\ell$ on $[-a,+\infty)$ that coincides with $W$ on $[-a,a-\eta]$. We also assume that $W_\ell$ is bounded (with all its derivatives) and that, on $[a-\eta,+\infty)$, $W_\ell\geq c>0$. We denote by $f_\ell$ the positive, $L^2$-normalized groundstate of the associated operator $\mathscr{L}_\ell=-h^2\partial_x^2+W_\ell$ and we denote by $\mu$ the corresponding eigenvalue. We consider a smooth cutoff function $0\leq \chi_\ell\leq 1$ supported in $[-a,a-\eta)$ and equal to $1$ on $[-a,a-2\eta)$.

Thanks to our decay estimate on $f_\ell$ and support considerations, we have
\[
(\mathscr{L}-\mu)(\chi_\ell f_\ell)=\mathcal{O}(e^{-(S_0-\varepsilon)/h})
\]
where $S_0=\int_{-a}^a\sqrt{W(s)}\mathrm{d}s$, and $\varepsilon >0$ can be made as small as we want by taking $\eta>0$ small enough. By symmetry, we also define the right-well potential $W_r$, and using self-evident notation, we have
\[(\mathscr{L}-\mu)(\chi_r f_r)=\mathcal{O}(e^{-(S_0-\varepsilon)/h})\,.\]
We also notice that
\[\langle\chi_\ell f_\ell,\chi_r f_r\rangle=\mathcal{O}(e^{-(S_0-\varepsilon)/h})\,.\]
Then, the spectral theorem shows that there are at least two eigenvalues near $\mu$ at a distance of order at most $\mathcal{O}(e^{-(S_0-\varepsilon)/h})$. In fact there are exactly two: we can show that the eigenfunctions of $\mathscr{L}$ are localized near $\pm a$ and that the low-lying spectra of $\mathscr{L}$ and $\mathscr{L}_\ell\oplus\mathscr{L}_r$ coincide modulo $\mathcal{O}(h^\infty)$. Due to the spectral asymptotics of $\mathscr{L}_{\ell/r}$ (see Theorem \ref{thm.linear.wells}), we deduce that 
\[\lambda_2(h)-\lambda_1(h)=\mathcal{O}(e^{-(S_0-\varepsilon)/h})\]
and that $\lambda_3(h)-\lambda_2(h)\geq ch^{\frac23}$.	

There remains to get the accurate asymptotics of the spectral gap. We denote by $\Pi$ the spectral projector onto ${\rm{Ker}}(\mathscr{L} - \lambda_1) \oplus {\rm{Ker}}(\mathscr{L} - \lambda_2)$. The above discussion and the spectral theorem imply that 
\begin{equation}\label{est.I-P} (\mathrm{Id}-\Pi)(\chi_\star f_\star)=\mathcal{O}(e^{-(S_0-\varepsilon)/h})\,.\end{equation}
We let 
\[g_\star=\Pi (\chi_\star f_\star)\,.\]
The vectors $(g_\ell,g_r)$ form a basis of ${\rm{Ran}} \Pi$. However, it is not orthonormal. Therefore, we consider the associated Gram matrix of the family $(g_\ell, g_r)$, namely the $2\times2$ matrix $G=(\langle g_j,g_k\rangle)$. Using \eqref{est.I-P} we can check that
\[G=\mathrm{Id}+T_1+\mathcal{O}(e^{-2(S_0-\varepsilon)/h})\,,\] with
\[ \quad  T_1=\begin{pmatrix}
0&\langle \chi_\ell f_\ell,\chi_r f_r\rangle\\
\langle \chi_r f_r,\chi_\ell f_\ell\rangle&0
\end{pmatrix}=\mathcal{O}(e^{-(S_0-\varepsilon)/h})\,,\]
which shows that the family $(g_\ell, g_r)$ is almost orthonormal. The family $G^{-\frac12}(g_\ell,g_r)$ is orthonormalized and the matrix of $\mathscr{L}$ in this basis is given by
\[L:=G^{-\frac12}\begin{pmatrix}
\langle\mathscr{L} g_\ell,g_\ell\rangle&\langle\mathscr{L} g_\ell,g_r\rangle\\
\langle\mathscr{L} g_r,g_\ell\rangle&\langle\mathscr{L} g_r,g_r\rangle\\
\end{pmatrix}G^{-\frac12}\,.\]
We have
\[\begin{pmatrix}
\langle\mathscr{L} g_\ell,g_\ell\rangle&\langle\mathscr{L} g_\ell,g_r\rangle\\
\langle\mathscr{L} g_r,g_\ell\rangle&\langle\mathscr{L} g_r,g_r\rangle\\
\end{pmatrix}=\mu\mathrm{Id}+T_2+\mathcal{O}(e^{-2(S_0-\varepsilon)/h})\,,\]
with
\[
T_2=\begin{pmatrix}
0&\langle\mathscr{L}(\chi_\ell f_\ell),(\chi_r f_r)\rangle\\
\langle\mathscr{L}(\chi_\ell f_\ell),(\chi_r f_r)\rangle&0
\end{pmatrix}\,.
\]
Therefore,
\[
L=\begin{pmatrix}
\mu&w\\
w&\mu
\end{pmatrix}+\mathcal{O}(e^{-2(S_0-\varepsilon)/h})\,,\quad w=\langle(\mathscr{L}-\mu)(\chi_\ell f_\ell),\chi_r f_r\rangle\,,
\]
and the spectral gap of this matrix is
\[
\lambda_2 - \lambda_1 = 2|w| + \mathcal{O}(e^{-2(S_0-\varepsilon)/h}).
\]
Then, we have, by support considerations and by using the definition of $f_\ell$,
\[
w=\langle(\mathscr{L}_\ell-\mu)(\chi_\ell f_\ell),\chi_r f_r\rangle=-h^2\langle[\partial_x^2,\chi_\ell]f_\ell,f_r\rangle=-h^2\langle2\chi'_\ell f'_\ell+\chi''_\ell f_\ell, f_r\rangle\,.
\]
By integration by parts, we deduce that
\begin{align*}
\langle2\chi'_\ell f'_\ell+\chi''_\ell f_\ell, f_r\rangle&=\langle2\chi'_\ell f'_\ell, f_r\rangle-\langle\chi'_\ell f'_\ell, f_r\rangle-\langle\chi'_\ell f_\ell, f'_r\rangle\\ &=\int_0^{a-\eta} \chi'_\ell(f'_\ell f_r-f_\ell f'_r)\mathrm{d}s\,,
\end{align*}
where $0$ is the middle point between the two wells.
Since, on $(0,a-\eta)$
\[
(f'_\ell f_r-f_\ell f'_r)'=f''_\ell f_r-f_\ell f''_r=h^{-2}(W_\ell-W_r)f_\ell f_r=0\,,
\]
a second integration by parts gives
\[
w=h^2(f'_\ell(0) f_r(0)-f_\ell(0) f'_r(0))=2h^2 f'_\ell(0)f_\ell(0)\,.
\]
Finally, we can use the pointwise estimate from Theorem \ref{thm.linear.wells}, for the single-well groundstate:
\begin{align*}
f_\ell(0) &\sim \frac{|W'(a)|^{\frac13}}{2\sqrt{\pi}|\mathrm{Ai}'(z_1)|}\frac{e^{-\frac{1}{h}S(0)}}{W(0)^{\frac14}h^{\frac16}}\\
f_\ell'(0)&\sim -\frac{\sqrt{W(0)}}{h} \frac{|W'(a)|^{\frac13}}{2\sqrt{\pi}|\mathrm{Ai}'(z_1)|}\frac{e^{-\frac{1}{h}S(0)}}{W(0)^{\frac14}h^{\frac16}}.
\end{align*}
We conclude that
\[
\lambda_2(h)-\lambda_1(h)\underset{h\to 0}{\sim} \frac{ h |W'(a)|^{2/3} }{\pi\mathrm{Ai}'(z_1)^2 h^{1/3}} e^{-S_h/h}\,.
\]
\end{proof}

\subsection{Adaptation to double $V$-shaped wells}
In this final subsection we discuss how to adapt the results of this article to $V$-shaped wells. Take $a>0$ and consider an even continuous non-negative electric potential $W$ that is piecewise smooth on $\R\setminus\{\pm a\}$ with exactly two minima at $\pm a$ such that $W(\pm a)=0$. Assume that $A=W'(a^+)>0$ and $-B=W'(a^-)<0$. Also assume that $W$ is confining, i.e $\liminf_\infty W >0$. We consider the corresponding operator 
\[
\mathscr{L}=-h^2\partial_x^2+W\,.
\]
The following operator plays an important role in the analysis:
\[
\mathscr{M}=-h^2\partial_x^2 - B x\ {1}_{\R_-}(x)+ A x\ {1}_{\R_+}(x)\,.
\]
It appears as a first order approximation of $\mathscr{L}$ near the well at $a$. This selfadjoint operator has compact resolvent and its eigenvalues and eigenfunctions can be expressed in terms of the Airy functions. By homogeneity, its groundstate eigenvalue $z(h)$ satisfies $z(h)=h^{\frac23}z(1)$.

\begin{proposition}\label{prop.application2}
There exists a $2/3$ homogeneous positive function $\kappa(A,B)$ such that
\[
\lambda_2(h)-\lambda_1(h)=\kappa(A,B) h^{\frac23}e^{-S_h/h}(1+o(1))\,,\quad S_h=\int_{-a}^a\sqrt{(W-h^{\frac23}z(1))_+}\,.
\]	
\end{proposition}

The $\kappa$ function is obtained by solving an implicit equation involving Airy functions. There does not seem to be a compact formula for it. 

\begin{proof}[Proof idea]
The proof follows the same lines as for Proposition \ref{prop.application1}. We start with a piecewise linear approximation near the bottom of a single well, which gives an $L^2$ approximation of the normalized eigenfunctions in terms of Airy functions (at leading order), and permits the computation of the relevant constants. Then a quasimode construction and Agmon estimates propagate the Airy approximation inside the forbidden region between the two wells. 
\end{proof}

\subsection*{Declaration of AI use}

AI tools were used solely to check the mathematical calculations and reasoning in the final version of the manuscript. They were not used to write the manuscript, develop proofs, or generate mathematical ideas.

	\bibliographystyle{abbrv}
	\bibliography{biblioMoRa}
	
	\end{document}